\documentclass{article}
\usepackage[english,russian,ukrainian]{babel}
\usepackage{latexsym,amssymb}
\usepackage[cp1251]{inputenc}
\usepackage{amsfonts,amssymb}
\usepackage{euscript}
\newcounter{theorem} 
\newcounter{lemma} 

\begin{document}

\vspace*{7mm}

\Large







\noindent {\bf  Vitalii Shpakivskyi}\\

\noindent {\bf Construction of an Infinite-Dimensional Family of Exact Solutions of  the Klein--Gordon Equation by the Hypercomplex Method}

\vspace{7mm}

The infinite-dimensional family of exact solutions of  the Klein--Gordon equation is constructed by the hypercomplex method.

\Large

\section{Introduction}

In quantized field theory, the Klein--Gordon (or Klein--Gordon--Fock) equation describing the spinor field created by charged particles with spin 1/2 is widely used. 
Several types of Klein--Gordon equation have been studied in the literature.
We will consider the Klein--Gordon equation of the form
 
 \begin{equation}\label{Kl-Gor}
\left(\Box_4+\mathfrak{m}\right)U(x_0,x_1,x_2,x_3):=\frac{\partial^2 U }{\partial
x_0^2}-\sum\limits_{j=1}^3\frac{\partial^2 U }{\partial
x_j^2}+\mathfrak{m} U=0,\qquad \mathfrak{m}\in\mathbb{C}.
\end{equation}

Equation (\ref{Kl-Gor}) is used to describe fast-moving particles that have mass (rest mass).
Many scientists have studied issues of constructing exact solutions of the Klein--Gordon equation. 
All works differ by types of Klein--Gordon equations, classes of solutions, and used methods.
Naturally, different methods give different exact solutions.

In papers \cite{Ahmadov,Burt,Ebaid,Wang,Yang,Changzheng,Vieira,Cote,Agom,Joseph,Cavalcanti,Yinghui,Olgar,Alvarez,
Dusunceli,Boldyreva,Sirendaoreji,Shakeri,Lakestani,Yao,Fu,Falaye,Yasuk,Alhaidari,Dib,Jia,Merad,Rezaei,Ikot}
  exact solutions of linear and nonlinear (generalized) Klein--Gordon equations were considered.

In this paper, by the hypercomplex method, we construct an infinite-dimensional family of exact solutions of equation (\ref{Kl-Gor}) of the form $$U_n(x_0,x_1,x_2,x_3)=P_n(x_0,x_1,x_2,x_3)\exp\Biggr(
f(x_0,x_1,x_2,x_3)\Biggr),$$ where $P_n$ is a complex polynomial of degree $n$, and $f(x_0,x_1,x_2,x_3)$ is a complex linear function of real variables $x_0,x_1,x_2,x_3$.

\section{The method}\label{Sect2}

This method originates from the works of P.\,W. Ketchum \cite{Ketchum-28,Ketchum-29-3-6,Ketchum-32-3-6}, 
and was independently developed in the studies of I.\,P. Mel'nichenko \cite{Mel'nichenko75,Melnichenko03}, 
S.\,A. Plaksa  \cite{Plaksa12} and ours.
The method development was finalized in our papers \cite{Shp-APAM,Shp-UMB-2018-3-4,Shp-ProcIGC-2018-3-5,Shp-ProcIAMM-2018-3-6}. 
This method is systematically described in monograph \cite[chapters 8, 11-13]{Plaksa-Shpakivskyi-2023}.

Let $\mathbb{A}$ be an $n$-dimensional commutative associative algebra
over the field of complex numbers $\mathbb{C}$ and let
$e_0,e_1,e_2,e_3\in \mathbb{A}$ be four linearly independent vectors over the real field $\mathbb{R}$. Let
$\zeta:=x_0e_0+x_1e_1+x_2e_2+x_3e_3$, $x_0,x_1,x_2,x_3\in\mathbb{R}$.  Let
 $E_4:=\{\zeta=x_0e_0+x_1e_1+x_2e_2+x_3e_3:\,\, x_0,x_1,x_2,x_3\in\mathbb{R}\}$ be the
linear span of vectors $e_0,e_1,e_2,e_3$ over the field
$\mathbb{R}$. Let $\Omega$ be a domain in $E_4$.

\vskip2mm

\textbf{Definition 1.} \textit{
We say that a continuous function
$\Phi:\Omega\rightarrow\mathbb{A}$ is \textit{monogenic}
in $\Omega$ if $\Phi$ is differentiable in the sense of
Gateaux in every point of $\Omega$, i.~e. if  for every
$\zeta\in\Omega$ there exists an element
$\Phi'(\zeta)\in\mathbb{A}$ such that
\begin{equation}\label{monogennaOZNA}\medskip
\lim\limits_{\varepsilon\rightarrow 0+0}
\frac{\Phi(\zeta+\varepsilon
h)-\Phi(\zeta)}{\varepsilon}= h\Phi'(\zeta)\quad\forall\,
h\in E_4.\medskip
\end{equation}
$\Phi'(\zeta)$ is the Gateaux derivative of the function
$\Phi$ in the point $\zeta$.}

In the algebra $\mathbb{A}$, we define the exponential function as the sum of the following
absolutely convergent series:
\begin{equation}\label{exp}
e^{\zeta}:=\sum\limits_{r=0}^\infty\frac{\zeta^r}{r!}\,.
\end{equation}
This function is evidently monogenic in $E_4$ and
 $$\frac{\partial^2 e^{\zeta}}{\partial x_j^2}=e^{\zeta}e_j^2\,,
\qquad j=0,1,2,3.$$

The consequence of equality (\ref{exp}) is equality 
\begin{equation}\label{exp-}
\left(\Box_4+\mathfrak{m}\right)e^\zeta=(e_0^2-e_1^2-e_2^2-e_3^2+\mathfrak{m})e^\zeta.
\end{equation}
By virtue of equality   (\ref{exp-}), the function $\Phi(\zeta)=e^\zeta$ satisfies equation (\ref{Kl-Gor}) if the vectors $e_0,e_1,e_2,e_3$ satisfy the relation
\begin{equation}\label{har-Kl-Gor}
e_0^2-e_1^2-e_2^2-e_3^2+\mathfrak{m}=0.
\end{equation}

Since the function $\Phi(\zeta)=e^{\zeta}$ satisfies equation (\ref{Kl-Gor}), all its complex-valued components in the decomposition with respect
to the basis of the algebra $\mathbb{A}$ are solutions of equation (\ref{Kl-Gor}).

Equation (\ref{har-Kl-Gor}) is called \textit{the characteristic equation} for equation (\ref{Kl-Gor}).

Thus, we have two problems. 

\textbf{Problem 1.} To find all solutions of characteristic equation (\ref{har-Kl-Gor}).

\textbf{Problem 2.}  To find the explicit form of the components of the monogenic function $\Phi(\zeta)=e^\zeta$.

\section{Solution to problem 1}\label{Sect3}

Instead of the algebra $\mathbb{A}$, we will consider \textit{a sequence} of algebras  $\{\mathbb{E}^n_\rho\}_{n=1}^\infty$ (see \cite[section 12.2]{Plaksa-Shpakivskyi-2023}). 
For any natural\, $n$, let\, $\mathbb{E}^n_\rho$\, be a commutative associative algebra
over the complex field with a basis\, $\{1,\rho,\rho^2,\ldots,\rho^{n-1}\}$, where\, $\rho^{n}=0$.

Let vectors\, $e_0,e_1,e_2,e_3\in\mathbb{E}^n_\rho$\,
have the following decompositions with respect to the basis of
the algebra\, $\mathbb{E}^n_\rho$:
\begin{equation}\label{e_0_e_1_e_2_e_3-k-ro}
 e_0:=\sum\limits_{r=0}^{n-1}k_r\,\rho^r\,,\quad e_1:=\sum\limits_{r=0}^{n-1}m_r\,\rho^r\,,\quad e_2:=\sum\limits_{r=0}^{n-1}g_r\,\rho^r\,,
 \quad e_3:=\sum\limits_{r=0}^{n-1}d_r\,\rho^r\,,
\end{equation}
where $ k_r,m_r,g_r,d_r\in\mathbb{C}$.

Using decomposition (\ref{e_0_e_1_e_2_e_3-k-ro}) and multiplication table of the algebra   $\mathbb{E}^n_\rho$, we have
$$
 e_0^2=\sum\limits_{r=0}^{n-1} B_r\rho^r\,,
$$
where
$$
B_0=k_0^2\,,\quad B_1=2k_0k_1\,,\quad B_2=k_1^2+2k_0k_2\,,
$$
and, in the general case,
\begin{equation}\label{3-6:5+KG}
B_r\equiv B_r(k_0,k_1,\ldots,k_r):=\left\{
\begin{array}{lrr}
&k_{r/2}^2+2\left(k_0k_r+k_1k_{r-1}+\cdots+k_{\frac{r}{2}-1}k_{\frac{r}{2}+1}\right)\\[1.5mm]
&\;\;\mbox{if}\;\;  r \;\;\mbox{is\, even},\\[1.5mm]
&2\left(k_0k_r+k_1k_{r-1}+\cdots+k_{\frac{r-1}{2}}k_{\frac{r+1}{2}} \right)\\[1.5mm]
&\;\;\mbox{if}\;\;  r \;\;\mbox{is\,odd}.
\end{array}
\right.
\end{equation}

Similarly, we have
$$
e_1^2=\sum\limits_{r=0}^{n-1} C_r\rho^r\,,
$$
where the coefficients\, $C_r$\, satisfies the relations
\begin{equation}\label{3-6:C_r=B_rKG}
C_r\equiv C_r(m_0,m_1,\ldots,m_r)=B_r(m_0,m_1,\ldots,m_r)
\end{equation}
for\, $r=0,1,\dots,n-1$.

Similarly,
$$
e_2^2=\sum\limits_{r=0}^{n-1} G_r\rho^r\,,
$$
where the coefficients\, $G_r$\, satisfies the relations
\begin{equation}\label{3-6:G_r=B_rKG}
G_r\equiv G_r(g_0,g_1,\ldots,g_r)=B_r(g_0,g_1,\ldots,g_r)
\end{equation}
for\, $r=0,1,\dots,n-1$.

Similarly,
$$
e_3^2=\sum\limits_{r=0}^{n-1} D_r\rho^r\,,
$$
where the coefficients\, $D_r$\, satisfies the relations
\begin{equation}\label{3-6:D_r=B_rKG}
D_r\equiv D_r(g_0,g_1,\ldots,g_r)=B_r(g_0,g_1,\ldots,g_r)
\end{equation}
for\, $r=0,1,\dots,n-1$.

Thus, characteristic equation (\ref{har-Kl-Gor}) is equivalent to the following system of equations
$$B_0-C_0-G_0-D_0+\mathfrak{m}=0,$$
\begin{equation}\label{har-syst-Kl-Gor}
B_r-C_r-G_r-D_r=0,\qquad r=1,2,\ldots,
\end{equation}
where $B_r,C_r,G_r,D_r$ are defined by equalities (\ref{3-6:5+KG}), (\ref{3-6:C_r=B_rKG}),
(\ref{3-6:G_r=B_rKG}), (\ref{3-6:D_r=B_rKG}),  respectively. 

According to Remark 13.1 \cite{Plaksa-Shpakivskyi-2023} in characteristic equation (\ref{har-Kl-Gor}), we can set the vectors\, $e_1, e_2, e_3$\, arbitrarily, i.e.,\,
$k_r$, $m_r$, and $g_r$  may be arbitrary complex numbers for  $r=0,1,\ldots, n-1$ for any natural\, $n$.
Then the vector $e_0$ can be found, and the coefficients $k_r$  can be found from system  (\ref{har-syst-Kl-Gor}).

We note that the first equation of (\ref{har-syst-Kl-Gor}) is quadratic with respect to $k_r$, and
for every fixed $r=1,2,\ldots$ equation (\ref{har-syst-Kl-Gor}) is linear with respect to $k_r$. 
Therefore, it is easy to find $k_r$  for each fixed $r=0,1,2,\ldots$.

So, characteristic equation  (\ref{har-Kl-Gor}) is solved on the sequence of algebras $\{\mathbb{E}^n_\rho\}_{n=1}^\infty$.

\section{Solution to problem 2}\label{Sect4}

Problem 2 is solved in our article \cite{Shpakiv-Zb-17-3-6} (see also monograph \cite[Section 13.3]{Plaksa-Shpakivskyi-2023}). 
We will briefly describe the main results of our case.

The next decomposition follows from decompositions (\ref{e_0_e_1_e_2_e_3-k-ro})
$$\zeta=x_0e_0+x_1e_1+x_2e_2+x_3e_3=\sum\limits_{r=0}^{n-1}(k_rx_0+m_rx_1+g_rx_2+d_rx_3)\,\rho^r=:\sum\limits_{r=0}^{n-1}\xi_r\,\rho^r,$$
where $\xi_0:=k_0x_0+m_0x_1+g_0x_2+d_0x_3$ and $\xi_r:=k_rx_0+m_rx_1+g_rx_2+d_rx_3$, $r=1,2,...,n-1$.
In this case, the decomposition of resolvent $(t-\zeta)^{-1}$ is of the form
\begin{equation}\label{rozkl-rezol-A_n^m}
(t-\zeta)^{-1}=\sum\limits_{r=0}^{n-1}A_r\,\rho^r\,\qquad\forall\,\zeta\in E_4 \quad \forall\,t\in\mathbb{C}:\,
t\neq \xi_0\,,
  \end{equation}
where the coefficients\, $A_r$\, are determined by
the following recurrence relations:
\begin{equation}\label{-rozkl-rezol-A_n^m}
A_0=\frac{1}{t-\xi_0}\,, \quad A_r=\frac{1}{t-\xi_0}(\xi_rA_0+\xi_{r-1}A_1+\cdots+\xi_1A_{r-1}),
\quad r=1,2,\ldots.
\end{equation}

In particular, we have the following initial functions:
$$A_0=\frac{1}{t-\xi_0}\,, \qquad A_1=\frac{\xi_1}{(t-\xi_0)^2},\qquad A_2=\frac{\xi_2}{(t-\xi_0)^2}+\frac{\xi_1^2}{(t-\xi_0)^3},\quad
$$
$$
A_3=\frac{\xi_3}{(t-\xi_0)^2}+\frac{2\xi_1\xi_2}{(t-\xi_0)^3}+\frac{\xi_1^3}{(t-\xi_0)^4},
$$
$$A_4=\frac{\xi_4}{(t-\xi_0)^2}+\frac{2\xi_1\xi_3+\xi_2^2}{(t-\xi_0)^3}+\frac{3\xi_1^2\xi_2}{(t-\xi_0)^4}+
\frac{\xi_1^4}{(t-\xi)^5}\,,
$$
$$A_5=\frac{\xi_5}{(t-\xi_0)^2}+\frac{2\xi_1\xi_4+2\xi_2\xi_3}{(t-\xi_0)^3}+\frac{3\xi_1^2\xi_3+3\xi_1\xi_2^2}{(t-\xi_0)^4}+
\frac{4\xi_1^3\xi_2}{(t-\xi_0)^5}+\frac{\xi_1^5}{(t-\xi_0)^6}\,,
$$
$$A_6=\frac{\xi_6}{(t-\xi_0)^2}+\frac{\xi_3^2+2\xi_1\xi_5+2\xi_2\xi_4}{(t-\xi_0)^3}+
\frac{\xi_2^3+6\xi_1\xi_2\xi_3+3\xi_1^2\xi_4}{(t-\xi_0)^4}$$
$$+\frac{4\xi_1^3\xi_3+6\xi_1^2\xi_2^2}{(t-\xi_0)^5}+\frac{5\xi_1^4\xi_2}{(t-\xi_0)^6}+\frac{\xi_1^6}{(t-\xi_0)^7}\,,
$$
etc.

Note that the definition of the
function\, $e^{\zeta}$\, in form (\ref{exp}) is equivalent to its determination
as the principal extension of the holomorphic function $e\,^t$
of a complex variable $t$ into the sequence $\mathbb{E}^n_\rho$:
\begin{equation}\label{3-6:exp}
e^{\zeta}=\frac{1}{2\pi i}\int\limits_\gamma e\,^t(t-\zeta)^{-1}\,dt,
 \end{equation}
where\, $\gamma$\, is an arbitrary closed rectifiable
Jordan curve in the complex plane that embraces the point\,\,
 $\xi_0=k_0x_0+m_0x_1+g_0x_2+d_0x_3$, see Subsection 3.1.2 in \cite{Plaksa-Shpakivskyi-2023} 
  or the monograph by E.~Hille and R.~Phillips
\cite[p. 165]{Hil_Filips-3-6}.

Therefore, as shown in \cite{Shp-ProcIAMM-2018-3-6} (see also subsections 13.2.2, 13.3 in \cite{Plaksa-Shpakivskyi-2023}), 
we have that an infinite-dimensional family of solutions of equation
(\ref{Kl-Gor}) on the space $\mathbb{R}^4$ is of the following form:
\begin{equation}\label{sim-2}
\biggl\{\frac{1}{2\pi i}\int\limits_\Gamma
e^tA_r\,dt\biggr\}_{r=0}^\infty\,,
\end{equation}
 where  $\Gamma$ is a closed Jordan rectifiable curve on the complex plane that surrounds the point
 $\xi_0$, and\, $A_r$\, are defined by the recurrence formulas (\ref{-rozkl-rezol-A_n^m}), where $r$ can increase to infinity.

Let us show a quick way to calculate integral solutions (\ref{sim-2}). To do this, we will introduce some definitions. 

Let $\varphi(t-\xi_0,\xi_1,\ldots,\xi_r)$ be an arbitrary complex function of $(r+1)$ complex variables. 
We define a linear operator $P$ associating with each function $\varphi$ a function of $r$ variables by  the rule
$$P\varphi\big((t-\xi_0)^s,\xi_1,\ldots,\xi_r\big)=\varphi\Big((s-1)!\,,
\xi_1,\ldots,\xi_r\Big)\qquad \forall\,s\in\{2,3,\ldots,r+1\}.
$$
For example, $$P\left(\frac{\xi_3}{(t-\xi_0)^2}+2\frac{\xi_1\xi_2}{(t-\xi_0)^3}+
\frac{\xi_1^3}{(t-\xi_0)^4}\right)=
\xi_3+\xi_1\xi_2+\frac{\xi_1^3}{3!}\,.$$

Now let us define the functions 
\begin{equation}\label{0}
 \widetilde{\mathbf{\mathcal{A}}}_0:=1,\quad \widetilde{\mathbf{\mathcal{A}}}_r(\xi_1,\xi_2,\ldots,\xi_r):=P\,A_r\Biggr((t-\xi_0)^s,\xi_1,\ldots,\xi_r\Biggr)
\end{equation}
$$ \forall\, s\in\{2,3,\ldots,r+1\}, \quad r=1,2,...,$$
where the functions $A_r$ are defined by (\ref{-rozkl-rezol-A_n^m}).
As noted in Remark 1 \cite{Shpakiv-Zb-17-3-6}, the functions $\widetilde{\mathbf{\mathcal{A}}}_r$ are related by the recurrence relations
\begin{equation}\label{0+}
\widetilde{\mathbf{\mathcal{A}}}_r=\xi_r\widetilde{\mathbf{\mathcal{A}}}_0+\frac{r-1}{r}\,\xi_{r-1}\widetilde{\mathbf{\mathcal{A}}}_1+
\frac{r-2}{r}\,\xi_{r-2}\widetilde{\mathbf{\mathcal{A}}}_2+\cdots+\frac{1}{r}\,
\xi_{1}\widetilde{\mathbf{\mathcal{A}}}_{r-1}\,.
\end{equation}

Then the family of solutions (\ref{sim-2}) can be written as (see Lemma 2 \cite{Shpakiv-Zb-17-3-6}):
\begin{equation}\label{sim-2+}
U_r(x_0,x_1,x_2,x_3):=\biggl\{e^{\xi_0}\,\widetilde{\mathbf{\mathcal{A}}}_r\biggr\}_{r=0}^\infty,
\end{equation}
where the polynomials $\widetilde{\mathbf{\mathcal{A}}}_r$ are defined by recurrence relations (\ref{0+}) (or (\ref{0})).

Let us write the first few solutions from the family (\ref{sim-2+})  (or (\ref{sim-2})).
We get
$$
U_0=e^{\xi_0},\qquad U_1=\xi_1e^{\xi_0}, \qquad
U_2=\left(\xi_2+\frac{\xi_1^2}{2!}\right)e^{\xi_0},
$$
$$U_3=\left(\xi_3+\xi_1\xi_2+\frac{\xi_1^3}{3!}\right)e^{\xi_0},$$
$$U_4=\biggr(\xi_4+\frac{2\xi_1\xi_3+\xi_2
^2}{2!}+\frac{3\xi_1^2\xi_2}{3!}+\frac{\xi_1^4}{4!}\biggr)e^{\xi_0},
$$
$$U_5=\biggr(\xi_5+\xi_1\xi_4+\xi_2\xi_3+\frac{\xi_1\xi_2^2+
\xi_1^2\xi_3}{2!}+\frac{\xi_1^3\xi_2}{3!}+\frac{\xi_1^5}{5!}\biggr)e^{\xi_0},
$$
etc.

Calculating the subsequent functions $\widetilde{\mathbf{\mathcal{A}}}_r$ by the recurrence formula
(\ref{0+}) (or (\ref{0})), we can write out the  infinite-dimensional family of solutions
 of equation (\ref{Kl-Gor}) that represent a family of form (\ref{sim-2+})  (or \ref{sim-2}).

Summarizing Sections \ref{Sect3} and \ref{Sect4}, we obtain the following theorem.

\vskip2mm

\textbf{Theorem.} \textit{In the space $\mathbb{R}^4$, equation (\ref{Kl-Gor}) has an infinite-dimensional family of exact solutions of the  form
(\ref{sim-2+})  (or (\ref{sim-2})),  
where the polynomials $\widetilde{\mathbf{\mathcal{A}}}_r$ are defined by recurrence relations (\ref{0+}) (or (\ref{0})), where  $r$ 
 can increase to infinity.  Moreover,
in coefficients $\widetilde{\mathbf{\mathcal{A}}}_r$ the variables $\xi_0$, $\xi_r$ are of the form $\xi_0 =k_0x_0+m_0x_1+g_0x_2+d_0x_3$, 
 $\xi_r=k_rx_0+m_rx_1+g_rx_2+d_rx_3$, $r=1,2,\ldots$, where $m_r,g_r,d_r$ for $r=0,1,2,\ldots$, are arbitrary
complex numbers, and $k_r$ for each fixed $r=0,1,2,\ldots$ are determined by equality  (\ref{har-syst-Kl-Gor}).
}

\textbf{Sketch of the proof.}\vskip2mm

\noindent 1. As shown at the beginning of Section \ref{Sect2}, when condition  (\ref{har-Kl-Gor}) is satisfied, the function $\Phi(\zeta)=e^\zeta$ 
 is a solution of equation  (\ref{Kl-Gor}).
Since equation (\ref{Kl-Gor}) is linear, all  complex-valued components of the function $\Phi(\zeta)=e^\zeta$ are also solutions of (\ref{Kl-Gor}).\vskip2mm

\noindent 2. In Section \ref{Sect3} it is proved that the characteristic equation  (\ref{har-Kl-Gor}) is equivalent to system (\ref{har-syst-Kl-Gor}).
 That is, the vectors $e_1,e_2,e_3$ in the variable $\zeta=x_0e_0+x_1e_1+x_2e_2+x_3e_3$ are arbitrary, and, moreover,
  the vector $e_0$ is determined by equalities (\ref{har-syst-Kl-Gor}).\vskip2mm

\noindent 3.
In the paper \cite{Shp-APAM}, it is proved that every monogenic function $\Phi$ with values 
in the algebra $\mathbb{E}_\rho^n$ has the form
\begin{equation}\label{mon}
\Phi(\zeta)=\sum\limits_{r=0}^{n-1}\rho^r\,\frac{1}{2\pi
i}\int\limits_{\Gamma} F_r(t)(t-\zeta)^{-1}\,dt,
\end{equation}
where $F_r$ are holomorphic functions of a complex variable and $\Gamma$ is a
closed Jordan rectifiable curve in $D$ which surrounds the point
$\xi_0$. In particular, the exponent $\Phi(\zeta)=e^\zeta$  has the form (\ref{3-6:exp}).

\vskip2mm

\noindent 4. In the paper \cite{Shp-ProcIAMM-2018-3-6} (see also \cite{Shpakiv-Zb-17-3-6}), it is shown that the expansion of 
the resolvent $(t-\zeta)^{-1}$ in the algebra $\mathbb{E}_\rho^n$  
has the form (\ref{rozkl-rezol-A_n^m}) with coefficients $A_r$ of the form (\ref{-rozkl-rezol-A_n^m}).

\vskip2mm

\noindent 5.
Substituting resolvent (\ref{rozkl-rezol-A_n^m}) into formula (\ref{3-6:exp}), we obtain a finite family of solutions of equation (\ref{Kl-Gor}):
\begin{equation}\label{sem-K-G}
\biggl\{\frac{1}{2\pi i}\int\limits_\Gamma
e^tA_r\,dt\biggr\}_{r=0}^n \qquad {\rm or}\qquad \biggl\{e^{\xi_0}\,\widetilde{\mathbf{\mathcal{A}}}_r\biggr\}_{r=0}^n.
 \end{equation}

\vskip2mm
\noindent 6. Finally, in the paper \cite{Shp-ProcIGC-2018-3-5}, it is substantiated that in a finite set of solutions (\ref{sem-K-G})  $n$\, can increase to infinity.

\section{Several first solutions from family  (\ref{sim-2+})}

Next we will write some first solutions from family (\ref{sim-2+}).

$\bullet$ \,\, For $r=0$ equation (\ref{har-syst-Kl-Gor}) has the form
\begin{equation}\label{g_0-K-G}
k_0^2-m_0^2-g_0^2-d_0^2+\mathfrak{m}=0 \qquad {\rm or} \qquad 
k_0=\pm  \sqrt{\mathfrak{m}+m_0^2+g_0^2+d_0^2}\,.
\end{equation}
Therefore, $\xi_0=\pm  \sqrt{\mathfrak{m}+m_0^2+g_0^2+d_0^2}\,\,x_0+m_0x_1+g_0x_2+d_0x_3$, where $m_0,g_0,d_0$ are arbitrary
complex numbers.

Thus, the \textbf{first exact solution} of equation (\ref{Kl-Gor})  from   family (\ref{sim-2+}) is of the form
$$
U_0(x_0,x_1,x_2,x_3)=e^{\xi_0}=\exp\left(\pm  \sqrt{\mathfrak{m}+m_0^2+g_0^2+d_0^2}\,\,x_0+m_0x_1+g_0x_2+d_0x_3\right),
$$
where $m_0,g_0,d_0$ are arbitrary
complex numbers.

$\bullet$ \,\,  For $r=1$ equation (\ref{har-syst-Kl-Gor}) has the following form
\begin{equation}\label{g_1-0-K-G}
2k_0k_1-2m_0m_1-2g_0g_1-2d_0d_1=0,
\end{equation}
i.e., we have the linear equation with respect to $k_1$.
Substituting expression (\ref{g_0-K-G}) into (\ref{g_1-0-K-G}), we obtain
$$
k_1=\pm \frac{m_0m_1+g_0g_1+d_0d_1}{ \sqrt{\mathfrak{m}+m_0^2+g_0^2+d_0^2} }\,.
$$
Therefore, 
$$
\xi_1=\pm \frac{m_0m_1+g_0g_1+d_0d_1}{ \sqrt{\mathfrak{m}+m_0^2+g_0^2+d_0^2} }\, x_0+m_1x_1+g_1x_2+d_1x_3\,,
$$
where $m_0,g_0,d_0,m_1,g_1,d_1$ are arbitrary
complex numbers.

Thus, the \textbf{second exact solution} of equation (\ref{Kl-Gor})  from   family (\ref{sim-2+})  is of the form
$$
U_1(x_0,x_1,x_2,x_3)=\xi_1e^{\xi_0}=\left[\pm \frac{m_0m_1+g_0g_1+d_0d_1}{ \sqrt{\mathfrak{m}+m_0^2+g_0^2+d_0^2} }\, x_0+m_1x_1+g_1x_2+d_1x_3\right]
$$
$$
\times\exp\left(\pm  \sqrt{\mathfrak{m}+m_0^2+g_0^2+d_0^2}\,\,x_0+m_0x_1+g_0x_2+d_0x_3\right),
$$
where $m_0,g_0,d_0,m_1,g_1,d_1$ are arbitrary
complex numbers.

$\bullet$ \,\,  
For $r=2$ equation (\ref{har-syst-Kl-Gor}) has the following form
\begin{equation}\label{g_2-0-K-G}
k_1^2+2k_0k_2-m_1^2-2m_0m_2-g_1^2-2g_0g_2-d_1^2-2d_0d_2=0
\end{equation}
i.e., we have the linear equation with respect to $k_2$.
From (\ref{g_2-0-K-G}) we obtain
$$
k_2=\pm\frac{1}{2\sqrt{\mathfrak{m}+m_0^2+g_0^2+d_0^2}}\Biggr[m_1^2+2m_0m_2+g_1^2+2g_0g_2+d_1^2+2d_0d_2
$$
$$
-\frac{(m_0m_1+g_0g_1+d_0d_1)^2}{\mathfrak{m}+m_0^2+g_0^2+d_0^2}\,\Biggr],
$$
where $m_0,g_0,d_0,m_1,g_1,d_1,m_2,g_2,d_2$ are arbitrary
complex numbers.

Therefore, 
$$
\xi_2=\pm\frac{1}{2\sqrt{\mathfrak{m}+m_0^2+g_0^2+d_0^2}}\Biggr[m_1^2+2m_0m_2+g_1^2+2g_0g_2+d_1^2+2d_0d_2
$$
$$
-\frac{(m_0m_1+g_0g_1+d_0d_1)^2}{\mathfrak{m}+m_0^2+g_0^2+d_0^2}\,\Biggr]\,x_0+m_2x_1+g_2x_2+d_2x_3\,,
$$
where $m_0,g_0,d_0,m_1,g_1,d_1,m_2,g_2,d_2$ are arbitrary
complex numbers.

Thus, the \textbf{third exact solution} of equation (\ref{Kl-Gor})  from   family (\ref{sim-2+})  is of the form
$$
U_2(x_0,x_1,x_2,x_3)=\left(\xi_2+\frac{\xi_1^2}{2!}\right)e^{\xi_0}=
$$
$$
=\left[\pm\frac{1}{2\sqrt{\mathfrak{m}+m_0^2+g_0^2+d_0^2}}\Biggr(m_1^2+2m_0m_2+g_1^2+2g_0g_2+d_1^2+2d_0d_2\right.
$$

$$
\left.-\frac{(m_0m_1+g_0g_1+d_0d_1)^2}{\mathfrak{m}+m_0^2+g_0^2+d_0^2}\,\right)\,x_0+m_2x_1+g_2x_2+d_2x_3
$$

$$
\left.+\frac{1}{2!}\left( \pm\frac{m_0m_1+g_0g_1+d_0d_1}{\sqrt{\mathfrak{m}+m_0^2+g_0^2+d_0^2}}\, x_0+m_1x_1+g_1x_2+d_1x_3\right)^2\right]
$$

$$
\times\exp\left( \pm \sqrt{\mathfrak{m}+m_0^2+g_0^2+d_0^2}\,\,x_0+m_0x_1+g_0x_2+d_0x_3\right),
$$
where $m_0,g_0,d_0,m_1,g_1,d_1,m_2,g_2,d_2$ are arbitrary
complex numbers.

Continuing in this way, we can write out an infinite number of exact solutions of equation  (\ref{Kl-Gor}).

\vskip2mm
\textbf{The algorithm.}
\begin{enumerate}
 \item 
   \begin{itemize}
          \item We set $r=0$;
            \item find $k_0$ and write out $\xi_0=k_0x_0+m_0x_1+g_0x_2+d_0x_3$ from the quadratic equation (\ref{g_0-K-G}).
           Here $m_0,g_0,d_0$  are arbitrary
complex numbers;            
  \item get a solution $U_0$ by substituting $\xi_0$ in (\ref{sim-2+}) under the condition $\widetilde{\mathbf{\mathcal{A}}}_0=1$.
        \end{itemize}
        
    \item 
    
    \begin{itemize}
          \item We put $r=1$;
            \item having $k_0$ from the linear equation (\ref{g_0-K-G}) we find $k_1$ and we write out $\xi_1=k_1x_0+m_1x_1+g_1x_2+d_1x_3$. 
            Here $m_0,g_0,d_0,m_1,g_1,d_1$ are arbitrary
complex numbers;            
  \item  $\xi_1$ we substitute in $\widetilde{\mathbf{\mathcal{A}}}_1$ from (\ref{0+});
  \item $\widetilde{\mathbf{\mathcal{A}}}_1$ we substitute in  (\ref{sim-2+}) and we obtain the solution $U_1$.
        \end{itemize}
        
        \item 
    
    \begin{itemize}
          \item We put $r=2$;
            \item having $k_0,k_1$ from the linear equation (\ref{g_0-K-G}) we find $k_2$ and we write out $\xi_2=k_2x_0+m_2x_1+g_2x_2+d_2x_3$.
             Here $m_0,g_0,d_0,m_1,g_1,d_1,m_2,g_2,d_2$ are arbitrary
complex numbers;
  \item  $\xi_1$ and $\xi_2$ we substitute in $\widetilde{\mathbf{\mathcal{A}}}_2$ from (\ref{0+});
  \item $\widetilde{\mathbf{\mathcal{A}}}_2$ we substitute in  (\ref{sim-2+}) and we obtain the solution $U_2$.
        \end{itemize}
              
       \item   etc.
        
\end{enumerate}

\section{Final remarks}

\textbf{1.} 
It is obvious that the proposed method works for Klein--Gordon equations of lower dimensions. For example, solutions of the three-dimensional equation
 \begin{equation}\label{Kl-Gor-3}
\left(\Box_3+\mathfrak{m}\right)U(x_0,x_1,x_2):=\frac{\partial^2 U }{\partial
x_0^2}-\frac{\partial^2 U }{\partial
x_1^2}-\frac{\partial^2 U }{\partial
x_2^2}+\mathfrak{m} U=0
\end{equation}
also have the form  (\ref{sim-2+})  (or (\ref{sim-2})). 
In this case, the solutions of equation (\ref{Kl-Gor-3}) can be obtained from the solutions of equation (\ref{Kl-Gor}) by setting $d_r=0$, $r=0,1,2\ldots$.
For example, for  equation (\ref{Kl-Gor-3})
$$
U_0(x_0,x_1,x_2)=e^{\xi_0}=\exp\left(\pm  \sqrt{\mathfrak{m}+m_0^2+g_0^2}\,\,x_0+m_0x_1+g_0x_2\right),
$$
where $m_0,g_0$ are arbitrary
complex numbers.
Other solutions are written in a similar way.

The same method is valid for a two-dimensional equation
$$
\left(\Box_2+\mathfrak{m}\right)U(x_0,x_1):=\frac{\partial^2 U }{\partial
x_0^2}-\frac{\partial^2 U }{\partial
x_1^2}+\mathfrak{m} U=0.
$$

\textbf{2.} 
Solutions of type (\ref{sim-2+})  (or (\ref{sim-2})) are called \textit{polynomial--exponential type solutions}. 
It is obvious that all solutions (\ref{sim-2+})  (or (\ref{sim-2})) are analytical.

 \subsection*{Acknowledgment}
This work was supported by a grant from the Simons Foundation
(1030291,V.S.Sh.).

\renewcommand{\refname}{References}

\vskip10mm
\noindent Vitalii Shpakivskyi

\noindent Institute of Mathematics of the\\
National Academy of Sciences of Ukraine, Kyiv

\noindent e-mail: shpakivskyi86@gmail.com

\end{document}